\documentclass[12pt]{article}
\usepackage{amsfonts, ifthen, theorem}
\usepackage{amsmath, amssymb, theorem}
\usepackage[dvips]{color}
\usepackage{mathrsfs}

\renewcommand\theenumi{\roman{enumi}}
\renewcommand\labelenumi{\textup{(\theenumi)}}

\definecolor{ownblue}{rgb}{0,0,0}
\definecolor{boxblue}{rgb}{0,0,0}
\definecolor{ownred}{rgb}{0,0,0}
\definecolor{owngreen}{rgb}{0,0,0}
\definecolor{owngreen2}{rgb}{0,0,0}

\DeclareFontFamily{OT1}{rsfs}{} \DeclareFontShape{OT1}{rsfs}{m}{n}{
<-7> rsfs5 <7-10> rsfs7 <10-> rsfs10}{}
\DeclareMathAlphabet\mathcurl{OT1}{rsfs}{m}{n}

\theoremheaderfont\scshape
\newtheorem{theorem}{Theorem}[section]
\newtheorem{lemma}[theorem]{Lemma}
\newtheorem{corollary}[theorem]{Corollary}
\newtheorem{definition}[theorem]{Definition}
\newtheorem{proposition}[theorem]{Proposition}
\newtheorem{example}[theorem]{Example}

\theorembodyfont{\rmfamily}
\newtheorem{remark}[theorem]{Remark}

\renewcommand\thefootnote{\fnsymbol{footnote}}

\newcommand\comment[1]{}
\newcommand{\deep}[1]{_{{}_{#1}}}
\newcommand\R{\mathbb R}
\newcommand\N{\mathbb N}
\newcommand\Z{\mathbb Z}
\newcommand\F{\mathcurl F}
\newcommand\E[2][\P]{\mathbb E_{#1}\left[ #2\right]}
\renewcommand\P{\mathbb P}
\newcommand\Q{\mathbb Q}
\newcommand\K{\mathcal K}
\renewcommand\d{\mathrm d}
\newcommand\qed{\hfill$\Box$}
\newcommand\ip[2]{\langle #1,#2 \rangle}
\newcommand\adm{\textup{adm}}
\newcommand\suprep{\overline\pi} 
\newcommand\radnik[1][\mathbb Q]{\frac{\d #1}{\d\mathbb P}}
\newcommand\normcl[2][\Q]{\overline{#2}^{L^1(#1)}}
\newcommand\ind{1\hspace{-2.5mm}1}
\newcommand\iind[1]{1\hspace{-#1mm}1} 
\newcommand\cone{\operatorname{cone}}

\newcommand\dju[4]{\bigcup_{#1}^{#2}\hspace{-#3mm\cdot}\hspace{#4mm}}
\newcommand\cdju[2]{\bigcup_{#1}^{#2}\hspace{-4.2mm\cdot}\hspace{2.2mm}}
\newcommand\ldju[2]{\bigcup_{#1}^{#2}\hspace{-9mm\cdot}\hspace{7mm}}
\newcommand\rsi{[\![}
\newcommand\lsi{]\!]}

\renewcommand\theequation{\thesection.\arabic{equation}}

\newcommand{\nohyphens}{\hyphenpenalty=10000\exhyphenpenalty=10000\relax}%
\makeatletter
\renewcommand{\@seccntformat}[1]{\csname the#1\endcsname.\hspace{1em}}%
\renewcommand\section{\@startsection{section}{1}{\z@}%
                {-3.5ex \@plus -1ex \@minus -.2ex}%
                {2.3ex \@plus.2ex}%
                {\setcounter{equation}0\bfseries\nohyphens}}%
\makeatother

\newenvironment{proof}[1][]{\noindent\textit{Proof#1.} }{\vskip\baselineskip}

\begin{document}
{\centering {\renewcommand\thefootnote{}\large\bfseries On random
measures, unordered sums and discontinuities of the first kind}
}
\par
\vspace{1em}
\begin{center}
{\scshape Frank Oertel}\\ \itshape Department of Mathematics\\
University College Cork
\end{center}

\renewcommand{\abstractname}{}
\begin{abstract}\renewcommand\thefootnote{}\footnotesize
By investigating in detail discontinuities of the first kind of
real-valued functions and the analysis of unordered sums, where the
summands are given by values of a positive real-valued function, we
develop a measure-theoretical framework which in particular allows
us to describe \textit{rigorously} the representation and meaning of
sums of jumps of type $\sum_{0 < s \leq t} \Phi \circ \vert \Delta
X_s \vert$, where $X : \Omega \times \R_+ \longrightarrow \R$ is a
stochastic process with regulated trajectories, $t \in \R_+$ and
$\Phi : \R_+ \longrightarrow \R_+$ is a strictly increasing function
which maps $0$ to $0$ (cf. Proposition \ref{prop:sum of jumps on R+
with invertible function}). Moreover, our approach enables a natural
extension of the jump measure of c\`{a}dl\`{a}g and adapted
processes to an integer-valued random measure of optional processes
with regulated trajectories which need not necessarily to be right-
or left-continuous (cf. Theorem \ref{thm:optional random measures}).
In doing so, we provide a detailed and constructive proof of the
fact that the set of all discontinuities of the first kind of a
given real-valued function on $\R$ is at most countable (cf. Lemma
\ref{lemma:right limits and left limits}, Theorem \ref{thm:at most
countably many jumps on compact intervals} and Theorem \ref{thm:at
most countably many jumps on R+}).

By using the powerful analysis of unordered sums, we hope that our
contributions fill an existing gap in the literature, since neither
a detailed proof of (the frequently used) Theorem \ref{thm:at most
countably many jumps on compact intervals} nor a precise definition
of sums of jumps seems to be available yet.
\end{abstract}

\section{Preliminaries and notations}
In this section, we introduce the basic notation and terminology
which we will throughout in this paper. To perpetuate the lucidity
of the main ideas, we only consider $\R$-valued functions and
$\R$-valued trajectories of stochastic processes, although a
transfer to the (finite) multi-dimensional case is easily possible.
Most of our notations and definitions including those ones
originating from the general theory of stochastic processes and
stochastic analysis are standard. We refer the reader to the
monographs \cite{ChWi90}, \cite{HeWaYa92}, \cite{Ja79},
\cite{JaSh03}, \cite{Kl98}, \cite{Me82} and \cite{Pr03}. Concerning
a basic introduction to the the powerful theory of unordered sums,
we recommend the monographs \cite{Di60} and \cite{PrMo91}. Since at
most countable unions of pairwise disjoint sets play an important
role in this paper, we use a symbolic abbreviation. For example, if
$A : = \bigcup_{n=1}^{\infty} A_n$, where $(A_n)_{n \in \N}$ is a
sequence of sets such that $A_i \cap A_j = \emptyset$ for all $i
\not= j$, we write shortly $A : = \dju{n=1}{\infty}{8.8}{7.5}A_n$.

Throughout this paper, $(\Omega , \mathcal{F}, {\bf{F}}, {\Bbb{P}})$
denotes a fixed probability space, together with a fixed filtration
${\bf{F}}$. Even if it is not explicitly emphasized, the filtration
${\bf{F}} = (\mathcal{F}_t)_{t \geq 0}$ always is supposed to
satisfy the usual conditions\footnote{$\mathcal{F}_0$ contains all
$\mathbb{P}$-null sets and ${\bf{F}}$ is right-continuous.}.
A real-valued (stochastic) process $X : \Omega \times \R_+
\longrightarrow \R$ (which may be identified with the family of
random variables $(X_t)_{t \geq 0}$, where $X_t(\omega) : =
X(\omega, t)$) is called \textit{adapted} (with respect to
${\bf{F}}$) if $X_t$ is $\F_t$-measurable for all $t \in \R_+$. $X$
is called \textit{right-continuous} (respectively
\textit{left-continuous}) if for all $\omega \in \Omega$ the
trajectory $X_{\bullet}(\omega) : \R_+ \longrightarrow \R, t \mapsto
X_t(\omega)$ is a right-continuous (respectively left-continuous)
real-valued function. If all trajectories of $X$ do have left-hand
limits (respectively right-hand limits) everywhere on $\R_+$, $X_{-}
= (X_{t -})_{t \geq 0}$ (respectively $X_{+} = (X_{t +})_{t \geq
0})$ denotes the \textit{left-hand} (respectively
\textit{right-hand}) \textit{limit process}, where $X_{0 -} : =
X_{0+}$ by convention. If all trajectories of $X$ do have left-hand
limits and right-hand limits everywhere on $\R_+$, the \textit{jump
process} $\Delta X = (\Delta X_t)_{t \geq 0}$ is well-defined on
$\Omega \times \R_+$. It is given by $\Delta X : = X_{+} - X_{-}$
(cf. also Section $2$). A right-continuous process whose
trajectories do have left limits everywhere on $\R_+$, is known as a
\textit{c\`{a}dl\`{a}g} process. If $X$ is $\mathcal{F} \otimes
\mathcal{B}({\Bbb{R}}_+)$-measurable, $X$ is said to be
\textit{measurable}. $X$ is said to be \textit{progressively
measurable} (or simply \textit{progressive}) if for each $t \geq 0$,
its restriction $X\vert_{\Omega \times [0, t]}$ is $\mathcal{F}_t
\otimes \mathcal{B}([0, t])$-measurable. Obviously, every
progressive process is measurable and (thanks to Fubini) adapted.

A random variable $T : \Omega \longrightarrow [0, \infty]$ is said
to be a \textit{stopping time} or \textit{optional time} (with
respect to ${\bf{F}}$) if for each $t \geq 0$, $\{T \leq t\} \in
{\mathscr{F}}_t$. Let $\mathcal{T}$ denote the set of all stopping
times, and let $S, T \in \mathcal{T}$ such that $S \leq T$. Then
$\rsi S, T \rsi : = \{ (\omega , t) \in \Omega \times \R_+ :
S(\omega) \leq t < T(\omega)\}$ is an example for a
\textit{stochastic interval}. Similarly, one defines the stochastic
intervals $\lsi S, T \lsi$, $\lsi S, T \rsi$ and $\rsi S, T \lsi$.
Note again that $\rsi T \lsi : = \rsi T, T \lsi =
\textup{Gr}(T)\vert_{\Omega \times \R_+}$ is simply the graph of the
stopping time $T : \Omega \longrightarrow [0, \infty]$ restricted to
$\Omega \times \R_+$. $\mathcal{O} = \sigma\big\{[\![T,\infty [\![
\hspace{1mm}: T \in \mathcal{T}\big\}$ denotes the \textit{optional
$\sigma$-field} which is generated by all c\`{a}dl\`{a}g adapted
processes. The \textit{predictable $\sigma$-field} $\mathcal{P}$ is
generated by all left-continuous adapted processes. An
$\mathcal{O}$- (respectively $\mathcal{P}$-) measurable process is
called \textit{optional} or \textit{well-measurable} (respectively
\textit{predictable}).
All optional or predictable processes are adapted. For the
convenience of the reader, we recall and summarise the precise
relation between those different types of processes in the following
\begin{theorem}\label{thm:POPA}
Let $(\Omega , \mathcal{F}, {\bf{F}}, {\Bbb{P}})$ be a filtered
probability space such that ${\bf{F}}$ satisfies the usual
conditions. Let $X$ be a stochastic process on $\Omega \times
{\Bbb{R}}_+$. Consider the following statements:
\begin{itemize}
\item[$($i$)$] $X$ is predictable;
\item[$($ii$)$] $X$ is optional;
\item[$($iii$)$] $X$ is progressive;
\item[$($iv$)$] $X$ is adapted.
\end{itemize}
Then the following implications hold:
\[
\textstyle{(i)} \Rightarrow \textstyle{(ii)} \Rightarrow
\textstyle{(iii)} \Rightarrow \textstyle{(iv)}.
\]
If $X$ is right-continuous, then the following implications hold:
\[
\textstyle{(i)} \Rightarrow \textstyle{(ii)} \iff \textstyle{(iii)}
\iff \textstyle{(iv)}.
\]
If $X$ is left-continuous, then all statements are equivalent.
\end{theorem}
\begin{proof}
The general chain of implications $\textstyle{(i)} \Rightarrow
\textstyle{(ii)} \Rightarrow \textstyle{(iii)} \Rightarrow
\textstyle{(iv)}$ is well-known (for a detailed discussion cf.
e.\,g. \cite{ChWi90}, Chapter 3). If $X$ is left-continuous and
adapted, then $X$ is predictable. Hence, in this case, all four
statements are equivalent. If $X$ is right-continuous and adapted,
then $X$ is optional (cf. \cite{ChWi90}, Remark following Theorem
3.4. and \cite{HeWaYa92}, Theorem 4.32). In particular, $X$ is
progressive. \qed
\end{proof}
By identifying processes that are almost everywhere identical, there
is no difference  between adapted \textit{measurable} processes,
optional processes, progressive processes and predictable processes
(cf. \cite{Ni99}). In particular, since every adapted
right-continuous process is optional, hence measurable, it is
therefore almost everywhere identical to a predictable process.

Let $A \subseteq \Omega \times {\Bbb{R}}_{+}$ and $\omega \in
\Omega$. Consider
\[
D_A(\omega) : = \inf \{t \in  {\Bbb{R}}_{+} : (\omega , t) \in A\}
\in [0, \infty]
\]
$D_A$ is said to be the \textit{d\'{e}but} of $A$. Recall that
$\inf(\emptyset) = + \infty$ by convention. $A$ is called a
\textit{progressive set} if $\ind_{A}$ is a progressively measurable
process. For a better understanding of the main ideas in the proof
of Theorem \ref{thm: jumps of optional processes and stopping
times}, we need the following non-trivial result (a detailed proof
of this statement can be found in e.\,g. \cite{Ba98} or
\cite{HeWaYa92}):
\begin{theorem}\label{thm:optional debut is stopping time}
Let $(\Omega , \mathcal{F}, {\bf{F}}, {\Bbb{P}})$ be a filtered
probability space such that ${\bf{F}}$ satisfies the usual
conditions. Let $A \subseteq \Omega \times {\Bbb{R}}_{+}$. If $A$ is
a progressive set, then $D_A$ is a stopping time.
\end{theorem}

\section{Discontinuities of the first kind}
In the following, let us denote by $I$ an arbitrary (bounded or
non-bounded) closed interval in $\R$, containing at least two
elements.
In other words, let $I$ be precisely one of the following sets:
\[
[a, b], [a, \infty), (-\infty, a], \R,
\]
where $a, b\in \R$, $a < b$. Let $f : I \longrightarrow {\Bbb{R}}$
be a real-valued function and $t \in I$ such that $(t, \infty) \cap
I \not= \emptyset$.\footnote{Any interior point of $I$ satisfies
that condition.}
Recall that the real value $f(t+)$ is the \textit{right-hand limit
of $f$} at $t$, if for every $\varepsilon > 0$ there exists a
$\delta > 0$ such that $\vert f(t+) - f(s) \vert < \varepsilon $
whenever $s \in (t, t + \delta)$.\footnote{Due to the choice of $t$
and the structure of $I$, we obviously may choose $\delta >0$
sufficiently small such that $(t, t + \delta) \subseteq I$.} Let $t
\in I$ such that $I \cap (-\infty, t) \not= \emptyset$. The real
value $f(t-)$ is said to be the \textit{left-hand limit of $f$} at
$b$, if for every $\varepsilon
> 0$ there exists a $\delta > 0$ such that $\vert f(t-) - f(s)
\vert < \varepsilon $ whenever $s \in (t - \delta, t)$. Let us
denote by
\[
L^+(f) : = \{t \in I : (t, \infty) \cap I \not= \emptyset \mbox{ and
}  f(t+) \textup{ exists}\}
\]
the set of all finite right-hand limits of $f$, and by
\[
L^-(f) : = \{t \in I : I \cap (-\infty, t) \not= \emptyset \mbox{
and } f(t-) \textup{ exists} \}
\]
the set of all finite left-hand limits of $f$. Let $t \in L(f) : =
L^+(f) \cap L^-(f)$. Then $\Delta f(t) : = f(t+) - f(t-) \in \R$
denotes the \textit{jump of $f$ at $t$}, leading to the well-defined
function $\Delta f : L(f) \longrightarrow \R$, the associated
\textit{function of jumps of $f$}. Let $\textup{int}(I)$ denote the
interior of $I$. An easy calculation shows that
\[
\textup{int}(I) = \{t \in I : (t, \infty) \cap I \not= \emptyset\}
\cap \{t \in I : I \cap (-\infty, t) \not= \emptyset\},
\]
and it follows that
\begin{equation}\label{eq:L(f)}
L(f) = \textup{int}(I) \cap \{t \in I : f(t-) \textup{ exists and }
f(t+) \textup{ exists}\}
\end{equation}
is a subset of $\textup{int}(I)$. The set $L(f)$ is known as the set
of \textit{discontinuities of the first kind} of $f$ or \textit{jump
points} of $f$. $I \setminus L(f)$ is called the set of
\textit{discontinuities of the second kind} of $f$ (cf.
\cite{Kl98}).

Fix an arbitrary $\varepsilon>0$ and consider the set $J(f;
\varepsilon)$ of all jumps of $f$ of size at least $\varepsilon$,
i.\,e.,
\[
J(f; \varepsilon) := \{t \in L(f) : \vert \Delta f(t)\vert \geq
\varepsilon\}.
\]
The set of \textit{all} jumps of the function $f$ is then given by
\[
J(f) : = \{t \in L(f) : \Delta f(t) \not= 0\} = \{t \in L(f) : \vert
\Delta f(t) \vert > 0\} = \bigcup\limits_{n \in \N}J(f;
\frac{1}{n}).
\]
Consider the function $\widetilde{f} : -I \longrightarrow
{\Bbb{R}}$, defined by $\widetilde{f}(s) : = f(-s)$. $\widetilde{f}$
simply describes the vertical reflection of $f$. Since the
right-hand limit of $f$ (respectively the left-hand limit of $f$) is
uniquely determined, vertical reflection of $f$ immediately implies
the following important\footnote{Note that the set $-I$ belongs to
the same class as the given interval $I$.}
\begin{proposition}\label{prop:reflection}
Let $f : I \longrightarrow {\Bbb{R}}$ be a real-valued function. Let
$\widetilde{f} : -I \longrightarrow {\Bbb{R}}$, defined by
$\widetilde{f}(s) : = f(-s)$ for all $s \in -I$. Then
\begin{itemize}
\item[$($i$)$]
$L^+(f) = -L^-(\widetilde{f})$, and $f(t+)= \widetilde{f}((-t)-)$
for all $t \in L^+(f)$;
\item[$($ii$)$]
$L^-(f) = -L^+(\widetilde{f})$, and $f(t-)= \widetilde{f}((-t)+)$
for all $t \in L^-(f)$.
\end{itemize} In particular,
$L(f) = -L(\widetilde{f})$ and
\[
J(f; \varepsilon) = - J(\widetilde{f}; \varepsilon)
\]
for all $\varepsilon >0$.
\end{proposition}
Clearly, there exists a direct link to the well-known and rich class
of regulated functions (cf. \cite{Di60}, 7.6. and \cite{Od04}). By
using our notation, recall that $f : I \longrightarrow \R$ is said
to be \textit{regulated on $I$} if and only if if $\textup{int}(I)
\subseteq \{t \in I : f(t-) \textup{ exists and } f(t+) \textup{
exists}\}$, if the left endpoint of $I$ belongs to $L^+(f)$, and if
the right endpoint of $I$ belongs to $L^-(f)$ (if the latter exist).
Consequently, due to \eqref{eq:L(f)}, we may state the following
\begin{remark}\label{remark: regulated fct}
Let $f : I \longrightarrow {\Bbb{R}}$ be a real-valued function.
Then the following statements are equivalent:
\begin{itemize}
\item[$($i$)$] f is regulated on $I$;
\item[$($ii$)$] $L(f) = \textup{int}(I)$, the left endpoint of $I$ belongs to
$L^+(f)$, and the right endpoint of $I$ belongs to $L^-(f)$ (if the
latter exist).
\end{itemize}
\end{remark}
Let $t \in \R$ and $(t_n)_{n\in \N}$ be an arbitrary sequence of
elements in a given non-empty subset $A$ of $\R$. If $\lim\limits_{n
\to \infty}t_n = t$ and $t < t_{n+1} < t_n$ (respectively $t_n <
t_{n+1} < t$) for all $n \in \N$, as usual, we make use of the
shorthand notation $t_n \downarrow t$ (respectively $t_n \uparrow
t$). Since compact intervals will play an important role later on,
the next statement is given for $I : = [a, b]$ only, where $a < b$.
However, as the proof clearly shows, our arguments are of local
nature, so that we actually may choose every interval $I$ of the
above type (including $\R_+$).
\begin{lemma}\label{lemma:right limits and left limits}
Let $a<b$, $f : [a, b] \longrightarrow {\Bbb{R}}$ be an arbitrary
real-valued function and $t \in [a, b]$.
\begin{itemize}
\item[$($i$)$] Let $(t_n)_{n \in \N} \subseteq L(f)$ such that
$t_n \downarrow t$. If $t \in L^+(f)$, then
\[
\lim\limits_{n \to \infty} f(t_n-) = f(t+) = \lim\limits_{n \to
\infty} f(t_n+).
\]
\item[$($ii$)$] Let $(t_n)_{n \in \N} \subseteq L(f)$ such that
$t_n \uparrow t$. If $t \in L^-(f)$, then
\[
\lim\limits_{n \to \infty} f(t_n-) = f(t-) = \lim\limits_{n \to
\infty} f(t_n+).
\]
\end{itemize}
In each of these cases, we have
\[
\lim\limits_{n \to \infty}\vert f(t_n+) - f({t_n-})\vert = 0.
\]
\end{lemma}
\begin{proof}
To verify $(i)$, let $t \in L^+(f)$ and $(t_n)_{n \in \N} \subseteq
L(f)$ such that $t_n \downarrow t$ and $n, m \in \N$ arbitrary. Put
$\tau_{mn} : = t_n - \xi_{mn}$, where $0 < \xi_{mn} : =
\frac{t_n-t_{n+1}}{2^m}$.
Then
\[ t < t_{n+1} < \tau_{mn} < t_n
\]
for all $m, n \in \N$, and $\tau_{mn} \uparrow t_n$ (as $m \to
\infty$) for all $n \in \N$.
Thus, using the definition of left-hand limits, we have
\begin{equation}\label{eqn:LHL} f({t_n}-) = \lim\limits_{m \to
\infty}f(\tau_{mn})
\end{equation}
for all $n \in \N$. Let $\varepsilon
> 0$. Since by assumption $t \in L^+(f)$, there exists a $\delta
> 0$ such that
\[
f\big((t, t + \delta)\big) \subseteq \big(f(t+) - \varepsilon, f(t+)
+ \varepsilon\big).
\]
Since $t_n \downarrow t$, it follows that $\lim\limits_{n \to
\infty}f(t_n) = f(t+)$ and that there exists $N_{\delta} \in \N$
such that $t_n - t < \delta$ for all $n \geq N_{\delta}$.
Consequently, $\tau_{mn} \in (t, t + \delta)$ for all $m \in \N$ and
$n \geq N_{\delta}$, implying that $\vert f(t+) - f(\tau_{mn})\vert
< \varepsilon$ for all for all $m, n \geq N_{\delta}$. In other
words, if $t \in L^+(f)$, then the double-sequence limit
$\lim\limits_{m,n \to \infty}f(\tau_{mn}) = f(t+)$ exists! Thanks to
a further epsilon-delta argument, we therefore obtain
\[ f(t+)
= \lim\limits_{n \to \infty}\big(\lim\limits_{m \to
\infty}f(\tau_{mn})\big) \stackrel{\eqref{eqn:LHL}}{=}
\lim\limits_{n \to \infty}f(t_n-).
\]
Now we use the same method to approach each $t_n$ decreasingly from
the right side. More precisely, let $m, n \in \N$ arbitrary and put
$\rho_{mn} : = t_n + \xi_{m,n-1}$, where $\xi_{m,0} := 0$
and $0 < \xi_{mn} : = \frac{t_n-t_{n+1}}{2^m}$.
Then
\[ t < t_n < \rho_{mn} < t_{n-1}
\]
for all $m \in \N, n \in \N \cap [2, \infty)$, and $\rho_{mn}
\downarrow t_n$ (as $m \to \infty$) for all $n \in \N$.
Thus, using the definition of right-hand limits, we have
\begin{equation}\label{eqn:RHL} f({t_n}+) = \lim\limits_{m \to
\infty}f(\rho_{mn})
\end{equation}
for all $n \in \N$. Again, since $t \in L^+(f)$, we obtain the
existence of a double-sequence limit, namely $\lim\limits_{m,n \to
\infty}f(\rho_{mn}) = f(t+)$. Hence,
\[ f(t+)
= \lim\limits_{n \to \infty}\big(\lim\limits_{m \to
\infty}f(\rho_{mn})\big) \stackrel{\eqref{eqn:RHL}}{=}
\lim\limits_{n \to \infty}f(t_n+).
\]
To complete the proof, we only have to consider the remaining case
$(ii)$. So, let $t \in L^-(f)$ and $(t_n)_{n \in \N} \subseteq L(f)$
such that $t_n \uparrow t$. Then $s_n \downarrow s$, where $s_n : =
-t_n$ and $s : = -t$. Consider the function $\widetilde{f} : [-b,
-a] \longrightarrow {\Bbb{R}}$, defined by $\widetilde{f}(s) : =
f(-s)$. Due to Proposition \ref{prop:reflection}, it follows that
$s_n = -t_n \in -L(f) = L(\widetilde{f})$ for all $n \in \N$ and $s
= -t \in -L^-(f) = L^+(\widetilde{f})$. Therefore, we precisely
obtain the situation of part $(i)$, but now related to the function
$\widetilde{f}$! Consequently,
\[
\lim\limits_{n \to \infty} \widetilde{f}(s_n-) = \widetilde{f}(s+) =
\lim\limits_{n \to \infty} \widetilde{f}(s_n+),
\]
and the claim follows by Proposition \ref{prop:reflection}. \qed
\end{proof}
If $f : \R_+ \longrightarrow \R$ is a regulated function, it follows
that $L^+(f) = \R_+$ and $L^{-}(f) = (0, \infty)$. Hence, we may
define $f_{+}(t) : = f(t+) \in \R$ for all $t \in \R_+$ and
$f_{-}(t) : = f(t-) \in \R$ for all $t \in (0, \infty)$, implying
the existence of well-defined functions $f_{+} : \R_+
\longrightarrow \R$ and $f_{-} : (0, \infty) \longrightarrow \R$. A
first immediate non-trivial implication of Lemma \ref{lemma:right
limits and left limits} is the following statement which will be
used in the proof of Lemma \ref{lemma:optional and regulated paths}.
\begin{corollary}\label{cor:f+ and f- for regulated f}
Let $f : \R_+ \longrightarrow \R$ be a regulated function. Then
$f_{+}$ is right-continuous on $\R_+$ and $f_{-}$ is left-continuous
on $(0, \infty)$.
\end{corollary}
\begin{theorem}\label{thm:at most countably many jumps on compact intervals}
Let $f : [a, b] \longrightarrow {\Bbb{R}}$ be an arbitrary
real-valued function, where $a<b$. Then
\begin{itemize} \item[$($i$)$] $J(f; \varepsilon)$ is finite for
all $\varepsilon >0$.
\item[$($ii$)$] $J(f)$ is at most countable.
\end{itemize}
\end{theorem}
\begin{proof}
Since $J(f) = \bigcup\limits_{n \in \N}J(f; \frac{1}{n})$, we only
have to prove $(i)$. Assume by contradiction that $J(f;
\varepsilon)$ is not finite. Due to the Bolzano-Weierstrass Theorem
the bounded and infinite set $J(f; \varepsilon)$ has at least one
accumulation point $t \in [a, b]$ (cf. e.\,g. \cite{Bu98}). Then
there exists a sequence $(t_n)_{n\in\N} \subseteq J(f; \varepsilon)$
such that $t_n \to t$ (as $n \to \infty$), $t_k \not= t_l$ for all
$k \not= l$, and $t_n \not= t$ for all $n \in \N$ (since $J(f;
\varepsilon)$ is not finite). We therefore can select a monotone
subsequence of $(t_n)_{n\in\N}$ which then also converges to $t$. To
avoid some cumbersome notation, WLOG, we may assume that the
original sequence $(t_n)_{n\in\N}$ is already the monotone one.
Consequently, we arrived exactly at either scenario $(i)$ or
scenario $(ii)$ of Lemma \ref{lemma:right limits and left limits}.
Since $t_n \in J(f; \varepsilon)$ for all $n \in \N$, we clearly
obtain a contradiction, and the claim follows. \qed
\end{proof}
The next result shows that at most countability of the jumps even
can be guaranteed for all real-valued functions which are defined on
the whole of $\R_+$ (respectively $\R$).
\begin{theorem}\label{thm:at most countably many jumps on R+}
Let $f : J \longrightarrow {\Bbb{R}}$ be an arbitrary real-valued
function, where $J\ \in \{\R_+, \R\}$. Then
\begin{itemize} \item[$($i$)$] $J(f; \varepsilon)$ is at most countable
for all $\varepsilon >0$.
\item[$($ii$)$] There exists a partition $\{D_k : k \in \N\}$
of $J(f)$ such that each $D_k$ is a finite subset of $J(f)$. In
particular, $J(f)$ is at most countable.
\end{itemize}
\end{theorem}
\begin{proof}
First, consider the case $J = \R_+$. Let $M : = \{t \in \R_+ : f(t-)
\textup{ exists and } f(t+) \textup{ exists}\}$. Since
$\textup{int}(\R_+) = (0, \infty) = \ldju{n=1}{\infty}(n-1, n]$,
representation \eqref{eq:L(f)} therefore implies that
\[
L(f) = \cdju{n=1}{\infty} \big((n-1, n] \cap M \big) =
\cdju{n=1}{\infty} \big( (n-1, n) \cap M \big) \cup (\N \cap M)
\stackrel{\eqref{eq:L(f)}}{=} \cdju{n=1}{\infty}L(f\vert_{[n-1, n]})
\cup (\N \cap M).
\]
Hence,
\begin{equation}\label{eq:jump size rep for functions on R+ II}
J(f) =  \cdju{n=1}{\infty}J(f\vert_{[n-1, n]}) \cup (\N \cap J(f))
\end{equation}
and
\begin{equation}\label{eq:jump size rep for functions on R+} J(f;
\varepsilon) =  \cdju{n=1}{\infty}J(f\vert_{[n-1, n]}; \varepsilon)
\cup (\N \cap J(f; \varepsilon))
\end{equation}
for all $\varepsilon > 0$. Thus, (i) follows by Theorem \ref{thm:at
most countably many jumps on compact intervals}. To prove (ii), fix
$n \in \N$ and consider $f_n : = f\vert_{[n-1, n]}$. Due to Theorem
\ref{thm:at most countably many jumps on compact intervals}, the set
$J(f_n; \frac{1}{m})$ is finite for each $m \in \N$. Since $J(f_n;
\frac{1}{m}) \subseteq J(f_n; \frac{1}{m+1})$ for all $m \in \N$, it
therefore follows that $J(f_n)$ can be written as an at most
countable union of disjoint \textit{finite} sets, namely as
\begin{equation}\label{eq:standard jump partition}
J(f_n) = \bigcup_{m = 1}^{\infty} J(f_n; \frac{1}{m}) = \dju{m =
1}{\infty}{4.5}{2.2} A_{m,n},
\end{equation}
where $A_{1, n} : = J(f_n; 1) = (\Delta f_n)^{-1}\big([1,
\infty)\big)$ and $A_{m+1, n} := J(f_n; \frac{1}{m+1}) \setminus
J(f_n; \frac{1}{m}) = (\Delta f_n)^{-1}\big([\frac{1}{m+1},
\frac{1}{m})\big)$ for all $m \in \N$. Hence, \eqref{eq:jump size
rep for functions on R+ II} implies that
\[
J(f) = \cdju{n=1}{\infty}\dju{m = 1}{\infty}{4.5}{2.2} A_{m,n} \cup
(\N \cap J(f)) = \cdju{k=1}{\infty} B_{k} \cup (\N \cap J(f)),
\]
where $\{B_k : k \in \N\} = \{A_{m,n} : (n,m) \in \N \times \N\}$.
Consequently,
\[
J(f) = \dju{l=1}{\infty}{3.8}{2.2} D_l,
\]
where $D_l : = B_l \cup (\{l\} \cap J(f))$ is a finite set for
all $l \in \N$.
\qed
\end{proof}
Since partitions of this type will play a fundamental role, we
introduce the following
\begin{definition}
Let $D$ be an at most countable subset of $\R$ which is not empty. A
partition $\{D_k : k \in \N\}$ of $D$ is called a finitely layered
partition of $D$ if $D = \dju{k = 1}{\infty}{8.8}{6.5} D_k$, where
$D_k$ is a finite subset of $D$ for all $k \in \N$.
\end{definition}
\comment{We continue with a recursive construction and an
observation which is easy to prove. Nevertheless, it shows
important applications in the general theory of stochastic
processes. Later on, by using finitely layered partitions, we will
recognise that it allows a recursive construction of stopping
times in a natural way (cf. proof of Theorem \ref{thm: jumps of
optional processes and stopping times}).
\begin{lemma}\label{lemma: countable sets}
Let $A$ be a finite subset of ${\Bbb{R}}$. Put
\[
s_1^A : = \min(A)
\]
and
\[
s_{n+1}^A := \min\big(A \cap (s_n^A, \infty)\big) = \min\big\{t \in
A : t > s_n^A \big\}
\]
for all $n \in \{1, 2, \ldots k-1\}$, where $k : =
\textup{card}(A)$. Then $s_n^A < s_{n+1}^A$ for all $n \in \{1, 2,
\ldots, k-1\}$, and
\[
A = \cdju{n = 1}k \{ s_n^A\}.
\]
\end{lemma}
}
\section{Unordered sums}
Using the previous results about the structure of the sets $L(f)$
and $J(f)$, we can introduce sum of jumps functions like e.\,g.
$\mathscr{P}(L(f)) \ni B \mapsto \sum_{s \in B} (\Delta f(s))^2 \in
[0, \infty]$ in a mathematically concise manner. Our aim is to
provide an exact description of such sums which is independent of
the choice of the partition of $J(f)$ (cf. Theorem
\ref{thm:unordered sum as measure} and Proposition \ref{prop:sum of
jumps on R+ with invertible function}). In particular, we will show
that finitely layered jump partitions provide a natural frame for
integer valued random measures
which are a \textit{special case} of such a (randomised) sum (cf.
Theorem \ref{thm:optional random measures}).

To this end, let $L$ be an arbitrary non-empty set and $h : L
\longrightarrow \R_+$ a \textit{positive} real-valued function.
Consider the set $\Bbb{F}(L) : = \{F : F \mbox { is a finite subset
of } L\}$. Clearly, $(\Bbb{F}(L), \subseteq \nolinebreak)$ is an
ordered set, and we may therefore consider the well-defined net $s_h
: \Bbb{F}(L) \longrightarrow \R_+$, defined by $s_h(F) : = \sum_{s
\in F} h(s)$, where $F \in \Bbb{F}(L)$. If the net $s_h$ converges
to a limit point $p \in \R_+$, $\sum_{s \in L} h(s) : = p$ is called
the \textit{unordered sum over $L$}. If the net $s_h$ converges,
$s_h$ is called \textit{summable}. Let us recall the following
\begin{theorem}\label{thm:summability}
Let $L$ be an arbitrary non-empty set and $h : L \longrightarrow
\R_+$ a positive real-valued function. Then the following statement
s are equivalent:
\begin{itemize}
\item[$($i$)$] $\sum_{s \in L} h(s)$ exists;
\item[$($ii$)$] The set $\{s_h(F) : F \in \Bbb{F}(L)\}$ is bounded in $\R_+$.
\end{itemize}
If the net $s_h$ converges, then $\sum_{s \in L} h(s) = \sup\{s_h(F)
: F \in \Bbb{F}(L) \}$.
\end{theorem}
Since we have to include the case that the net $s_h$ is not
convergent, Theorem \ref{thm:summability} justifies the following
natural extension of the unordered sum above:
\begin{definition}
Let $L$ be an arbitrary non-empty set and $h : L \longrightarrow
\R_+$ an arbitrary positive real-valued function. Define
\[
\sum_{s \in L} h(s) : = \sup\Big\{\sum_{s \in F} h(s) : F \in
\Bbb{F}(L) \Big\} 
\]
If $\emptyset \not= A \subseteq L$, put $\sum_{s \in A} h(s) : =
\sum_{s \in A} h\vert_A(s)$. Put $\sum_{s \in \emptyset} h(s) : =
0$.
\end{definition}
First note that in general, $\sum_{s \in L} h(s) \in [0, \infty]$
and that $\sum_{s \in E} h(s) \leq \sum_{s \in F} h(s)$ for all
subsets $E \subseteq F \subseteq L$. If $L = \{s_1, \ldots, s_n\}$
itself is a finite set, then obviously $\sum_{s \in L} h(s) =
\sum_{i=1}^{n}h(s_i) = \sum_{i=1}^{n}h(s_{\sigma(i)})$ for all
permutations $\sigma \in S_n$, which justifies the notation.
However, the following important fact, which we will use later on,
requires a proof.
\begin{lemma}\label{lemma:indicator function in sum}
Let $L$ be an arbitrary non-empty set and $h : L \longrightarrow
\R_+$ an arbitrary positive real-valued function. Let $A$ and $B$ be
arbitrary subsets of $L$. Then the following statements hold:
\begin{itemize}
\item[$($i$)$]
$\sum_{s \in L} h(s) \ind_A(s) < +\infty$ if and only if $\sum_{s
\in A} h(s) < +\infty$. Moreover,
\[
\sum_{s \in A} h(s) = \sum_{s \in L} h(s) \ind_A(s).
\]
\item[$($ii$)$]
$\sum_{s \in A} h(s) \ind_B(s) < +\infty$ if and only if $\sum_{s
\in A \cap B} h(s) < +\infty$. Moreover,
\[
\sum_{s \in A \cap B} h(s) =  \sum_{s \in A} h(s)\ind_{B}(s).
\]
\end{itemize}
These sums may be finite or infinite.
\end{lemma}
\begin{proof}
Since (ii) obviously follows by (i) (by applying (i) to the function
$h \ind_B$), we only have to prove (i). If $A = \emptyset$, nothing
is to prove. So, let $A \not= \emptyset$. Assume first that $\sum_{s
\in A} h(s) < +\infty$. Let $F$ be an arbitrary finite subset of
$L$. Since $F$ equals the disjoint union of the (finite) sets $A
\cap F$ and $(L\setminus A) \cap F$, standard associative and
commutative summation of finitely many numbers immediately gives
\[
\sum_{s \in F} h(s) \ind_A(s) = \sum_{s \in A \cap F} h(s) \ind_A(s)
= \sum_{s \in A \cap F} h(s).
\]
Since the finite subset $F$ of $L$ was arbitrarily chosen, it
therefore follows that
\[
\sum_{s \in L} h(s) \ind_A(s) \leq \sum_{s \in A} h(s) < +\infty.
\]
Now let  $\sum_{s \in L} h(s) \ind_A(s)$ be finite. Then, if $G$ is
an arbitrary finite subset of $A \subseteq L$, we obviously have
\[
\sum_{s \in G} h(s) = \sum_{s \in G} h(s) \ind_A(s) \leq \sum_{s \in
L} h(s) \ind_A(s) < +\infty,
\]
which proves the other inequality. Consequently, we have shown that
the equality holds if $\sum_{s \in A} h(s)$ is finite or if $\sum_{s
\in L} h(s) \ind_A(s)$ is finite. Hence, it must be true if these
sums are finite or if this is not the case. \qed
\end{proof}
\begin{proposition}\label{prop:linearity of unordered sums}
Let $L$ be an arbitrary non-empty set, $\alpha, \beta \geq 0$ and
$h, g : L \longrightarrow \R_+$ arbitrary positive real-valued
functions. Then $\sum_{s \in L} g(s) < +\infty$ and $\sum_{s \in L}
h(s) < +\infty$ if and only if $\sum_{s \in L} (\alpha g(s) + \beta
h(s)) < +\infty$. Moreover,
\[
\sum_{s \in L} (\alpha g(s) + \beta h(s)) = \alpha \sum_{s \in L}
g(s) + \beta \sum_{s \in L} h(s).
\]
These sums may be finite or infinite.
\end{proposition}
\begin{proof}
First, let $\sum_{s \in L} g(s) < +\infty$ and $\sum_{s \in L} h(s)
< +\infty$. Since addition is associative and commutative, the
equality obviously is true for every finite subset $F$ of $L$.
Consequently, we already obtain the inequality
\[
\sum_{s \in L} (\alpha g(s) + \beta h(s)) \leq \alpha \sum_{s \in L}
g(s) + \beta \sum_{s \in L} h(s) < +\infty.
\]
Now let $\sum_{s \in L} (\alpha g(s) + \beta h(s)) < +\infty$. Let
$E$ be an arbitrary finite subset of $L$. Then,
\[
\max\Bigg\{\alpha \sum_{s \in E} g(s), \beta \sum_{s \in E}
g(s)\Bigg\} \leq \sum_{s \in E} (\alpha g(s) + \beta h(s)) \leq
\sum_{s \in L} (\alpha g(s) + \beta h(s)) < +\infty,
\]
and it follows that both, $\Gamma : = \sum_{s \in L} g(s)$ and
$\Delta : = \sum_{s \in L} h(s)$ are finite. To prove the other
inequality, let $\varepsilon > 0$ be given. Then there exist finite
subsets $F$ and $G$ of $L$ such that
\[ \alpha \Gamma + \beta \Delta
< \sum_{s \in F} \alpha g(s) + \sum_{s \in G} \beta h(s) +
\varepsilon.
\]
Since we currently are working with summation of finitely many
elements only, we obviously may conclude that
\[
\alpha \Gamma + \beta \Delta - \varepsilon < \sum_{s \in F \cup G}
\alpha g(s) + \sum_{s \in F \cup G} \beta h(s) = \sum_{s \in F \cup
G} (\alpha g(s) + \beta h(s)).
\]
Since $F \cup G$ is a finite subset of $L$, we have arrived at the
other inequality. Hence, similarly as in the proof of Lemma
\ref{lemma:indicator function in sum}, we have shown that the
equality holds if $\sum_{s \in L} (\alpha g(s) + \beta h(s)) <
+\infty$ or if $\sum_{s \in L} g(s) < +\infty$ and $\sum_{s \in L}
h(s) < +\infty$. \qed
\end{proof}
\begin{corollary}\label{cor:disjoint sums}
Let $L$ be an arbitrary non-empty set and $h : L \longrightarrow
\R_+$ an arbitrary positive real-valued function. Let $C$ and $D$ be
arbitrary subsets of $L$. If $C \cap D = \emptyset$, then
\[
\sum_{s \in C \cup D} h(s) = \sum_{s \in C} h(s) + \sum_{s \in D}
h(s).
\]
These sums may be finite or infinite.
\end{corollary}
\begin{proof}
Since $C \cap D = \emptyset$, we have $\ind_{C \cup D} = \ind_{C} +
\ind_{D}$. Hence, by Lemma \ref{lemma:indicator function in sum} and
Proposition \ref{prop:linearity of unordered sums}, it therefore
follows that
\[
\sum_{s \in C \cup D} h(s) = \sum_{s \in L} h(s)\ind_{C \cup D}(s) =
\sum_{s \in L} (h(s)\ind_{C}(s) + h(s)\ind_{D}(s)) = \sum_{s \in C}
h(s) + \sum_{s \in D} h(s).
\]
\qed
\end{proof}
\begin{theorem}\label{thm:Fubini for unordered sums}
Let $L$ be an arbitrary non-empty set and $h : L \longrightarrow
\R_+$ a real-valued function. Let $\{D_n : n\in \N\}$ be an
arbitrary partition of a set $D \subseteq L$. Then the following
statements are equivalent:
\begin{itemize}
\item[$($i$)$] $\sum_{s \in D} h(s) < + \infty$;
\item[$($ii$)$] $\sum_{s \in D_n} h(s) < +\infty$ for all $n \in \N$ and
$\Big( \sum_{s \in D_n} h(s)\Big)_{n \in \N} \in l_1$.
\end{itemize}
Moreover,
\[
\sum_{s \in D} h(s) = \sum_{n = 1}^{\infty}\sum_{s \in D_n} h(s).
\]
These sums may be finite or infinite.
\end{theorem}
\begin{proof}
Nothing is to show if $D = \emptyset$. So, let $D \not= \emptyset$,
and assume first that (i) holds. Since $D_n \subseteq D$ for all $n
\in \N$, each finite subset of each $D_n$ is already a finite subset
of $D$, implying that $\sum_{s \in D_n} h(s) \leq \sum_{s \in D}
h(s) < +\infty$ for all $n \in \N$. Let $n \in \N$ arbitrary and
consider the set $C_n : = \ldju{k = 1}{n} D_k \subseteq D$. Due to
Corollary \ref{cor:disjoint sums}, we have
\[
0 \leq \sum_{k = 1}^{n}\sum_{s \in D_k} h(s) = \sum_{s \in C_n} h(s)
\leq \sum_{s \in D} h(s) < +\infty.
\]
Since $n \in \N$ was arbitrarily chosen, we may conclude that
\[
0 \leq \sum_{k = 1}^{\infty}\sum_{s \in D_k} h(s) \leq \sum_{s \in
D} h(s) < +\infty.
\]
Hence, $\Big( \sum_{s \in D_k} h(s)\Big)_{k \in \N} \in l_1$, and
statement (ii) follows. Now assume that (ii) holds. Then $0 \leq
\sum_{s \in D_n} h(s) < +\infty$ for all $n \in \N$ and $0 \leq
\sum_{n = 1}^{\infty}\sum_{s \in D_n} h(s) < + \infty$. Let F be an
arbitrary finite subset of $D$. Choose a sufficiently large number
$n \in \N$ such that $F \subseteq \ldju{k = 1}{n}D_k$, implying that
$F = \ldju{k = 1}{n}F_k$, where $F_k : = F \cap D_k$ for all $k \in
\{1, 2, \cdots , n\}$. Consequently, we have
\[
\sum_{s \in F} h(s) = \sum_{k = 1}^{n}\sum_{s \in F_k} h(s).
\]
Since each $F_k$ is a finite subset of $D_k$, assumption (ii)
further implies that
\[
\sum_{k = 1}^{n}\sum_{s \in F_k} h(s) \leq \sum_{k = 1}^{n}\sum_{s
\in D_k} h(s) \leq \sum_{k = 1}^{\infty}\sum_{s \in D_k} h(s) <
+\infty.
\]
Since the finite subset $F$ of $D$ was arbitrarily chosen, it
follows that statement (i) is true, and we have
\[
\sum_{s \in D} h(s) \leq  \sum_{n = 1}^{\infty}\sum_{s \in D_n} h(s)
< +\infty.
\]
Clearly, we have shown that the equality holds if the case (i) or
the case (ii) is given. Since (i) is equivalent to (ii), the
equality necessarily also must hold if one of the both unordered
sums is not finite. \qed
\end{proof}
Since $\N \times \N = \dju{m=1}{\infty}{9.8}{7.5}\ldju{n=1}{\infty}
\{(m, n)\} = \ldju{n=1}{\infty}\dju{m=1}{\infty}{9.8}{7.8} \{(m,
n)\}$, Theorem \ref{thm:Fubini for unordered sums} immediately
recovers a well-known result concerning the rearrangement of the
terms in a double series (cf. e.\,g. \cite{Bu98}):
\begin{corollary}\label{cor:double series}
Let $(a_{mn})_{(m, n) \in \N \times \N}$ be an arbitrary
double-sequence in $\R_+$. Then
\[
\sum_{m=1}^{\infty} \sum_{n=1}^{\infty} a_{mn} = \sum_{n \in \N
\times \N} a_{mn} = \sum_{n=1}^{\infty} \sum_{n=1}^{\infty} a_{mn}
\]
\end{corollary}
By using the language of measure theory, we have proven the
following important result:
\begin{theorem}\label{thm:unordered sum as measure}
Let $L$ be an arbitrary non-empty set and $h : L \longrightarrow
\R_+$ an arbitrary positive real-valued function. Then
\begin{eqnarray*}
\mu_h : \mathscr{P}(L) & \longrightarrow & [0,
\infty]\\
A & \mapsto & \sum_{s \in A} h(s),
\end{eqnarray*}
is a well-defined measure on the measurable space $(L,
\mathscr{P}(L))$.
\end{theorem}
\begin{remark}\label{measure theoretic version of indicator fct in sum}
Let $A$ and $B$ be arbitrary subsets of $L$. Then Lemma
\ref{lemma:indicator function in sum} implies that
\begin{equation}\label{eq:mu of A versus mu of L}
\mu_h(A) = \sum_{s \in L} h(s)\ind_A(s) =
\mu_{h\iind{2}_A}(L)
\end{equation}
and
\begin{equation}\label{eq:mu of B cap D}
\mu_h(A \cap B) = \sum_{s \in A} h(s)\ind_B(s) =
\mu_{h\iind{2}_B}(A).
\end{equation}
\end{remark}
Dependent on the choice of the function $h$, we recognise two
special and well-known cases:
\begin{itemize}
\item[$($i$)$] If $h(s) : = 1 = \ind_L(s)$ for all $s \in L$, then
$\mu_{\iind{2}_{L}}$ is precisely the {\it counting measure}.
\item[$($ii$)$] If $h : = \ind_{\{s_0\}}$, where $s_0 \in L$,
we obtain exactly the {\it Dirac measure} at $s_0$, since
\[
\mu_{\iind{2}_{\{s_0\}}}(A) \stackrel{\eqref{eq:mu of B cap D}}{=}
\mu_{\iind{2}_{L}}(\{s_0\} \cap A) \stackrel{\eqref{eq:mu of B cap
D}}{=} \mu_{\iind{2}_{A}}(\{s_0\}) = \ind_A(s_0) = \delta_{s_0}(A)
\]
for all $A \in \mathscr{P}(L)$.
\end{itemize}
\begin{corollary}\label{cor:countable sums}
Let $L$ be an arbitrary non-empty set and $h : L \longrightarrow
\R_+$ an arbitrary positive real-valued function.
\begin{itemize}
\item[$($i$)$] If $A \in \mathscr{P}(L)$ is finite, then
\[
\mu_h(A) = \sum_{\nu = 1}^{n}h(\nu),
\]
where $n = card(A)$.
\item[$($ii$)$] If $A \in \mathscr{P}(L)$ is countable, then
\[
\mu_h(A) = \sum_{n = 1}^{\infty}h(\varphi(n))
\]
for all bijective mappings $\varphi : \N \longrightarrow A$.
\end{itemize}
\end{corollary}
\begin{proof}
First note that $\mu_h(\{a\}) = h(a)$ for all $a \in A \subseteq L$.
Statement (i) now follows directly by Theorem \ref{thm:unordered sum
as measure}. To prove (ii), let $a_n : = \varphi(n)$, where $n \in
\N$ is arbitrary. Then, by Theorem \ref{thm:unordered sum as
measure} again, we have
\[
\mu_h(A) = \mu_h(\cdju{n=1}{\infty} \{a_n\}) = \sum_{n=1}^{\infty}
\mu_h(\{a_n\}) = \sum_{n=1}^{\infty} h(a_n),
\]
and the proof is finished. \qed
\end{proof}
We have developed all necessary tools which now allow us to give a
lucid and short proof of the following non-trivial result.
\begin{theorem}\label{thm:countably many positive values}
Let $L$ be an arbitrary non-empty set and $h : L \longrightarrow
\R_+$ an arbitrary positive real-valued function. Put $P : = \{s \in
L : h(s) > 0\}$. If $\sum_{s \in L}h(s) < +\infty$, then $P$ is at
most countable, and
\[
\sum_{s \in L}h(s) =  \sum_{s \in P}h(s) = \sum_{n =
1}^{\infty}h(\varphi(n))
\]
for all bijective mappings $\varphi : \N \longrightarrow P$.
\end{theorem}
\begin{proof}
By assumption, $\Sigma : = \sum_{s \in L}h(s) < +\infty$. Since $P =
\bigcup_{n = 1}^{\infty}P_n$, where $P_n : = \{s \in L : h(s) >
\frac{1}{n}\}$, we only have to show that each subset $P_n$ of $L$
is at most countable. We even show more and claim
that
\begin{equation}\label{eq:P_n is finite}
P_n \mbox{ is finite and consists of at most } \lfloor n\Sigma
\rfloor \mbox{ elements for all } n \in \N,
\end{equation}
where $\R \ni x \mapsto \lfloor x \rfloor : = \max\{m \in \Z : m
\leq x\}$ describes the assignment rule of the floor function. We
assume by contradiction that \eqref{eq:P_n is finite} is false. Then
there would exist $m \in \N$ and a \textit{finite} subset $G_m$ of
$P_m$ such that $\mbox{card}({G_m}) = \lfloor m\Sigma \rfloor + 1$.
But then, due to the definition of the floor function, we would have
\[
+\infty > \lfloor m\Sigma \rfloor + 1 > m\Sigma = \sum_{s \in L}m
h(s) \geq \sum_{s \in G_m}m h(s) > \mbox{card}(G_m)\cdot 1 = \lfloor
m\Sigma \rfloor + 1,
\]
which obviously is a contradiction. Hence, statement \eqref{eq:P_n
is finite} is true, implying that the set $P$ is at most countable.

Clearly, we have $\sum_{s \in L\setminus P}h(s) = 0$. Hence,
$\sum_{s \in L}h(s) = \sum_{s \in P}h(s) = \mu_h(P)$ (due to
Corollary \ref{cor:disjoint sums}), and Corollary \ref{cor:countable
sums} finishes the proof. \qed
\end{proof}
By linking Lemma \ref{lemma:indicator function in sum} and (the
proof of) Theorem \ref{thm:Fubini for unordered sums}, we can
characterise the finiteness of the measure $\mu_h$ in the following
way:
\begin{proposition}\label{prop:when is mu_h finite?}
Let $L$ be an arbitrary non-empty set and $h : L \longrightarrow
\R_+$ an arbitrary positive real-valued function. Then the following
statements are equivalent:
\begin{itemize}
\item[$($i$)$] $\mu_h : \mathscr{P}(L) \longrightarrow \R_+$ is a finite
measure.
\item[$($ii$)$] If $\{L_n : n \in \N\}$ is an arbitrary partition of
$L$ such that $\sum_{s \in L_n} h(s) < +\infty$ for all $n \in \N$
and $\Big(\sum_{s \in L_n} h(s)\Big)_{n \in \N} \in l_1$, then
\[
\mu_h(A) = \sum_{n = 1}^{\infty}\sum_{s \in L_n} h(s)\ind_{A}(s)
\]
for all $A \in \mathscr{P}(L)$.
\item[$($iii$)$] There exists a partition $\{C_l : l \in \N\}$ of
$L$ such that $\sum_{s \in C_l} h(s) < +\infty$ for all $l \in \N$,
$\Big(\sum_{s \in C_l} h(s)\Big)_{l \in \N} \in l_1$ and
\[
\mu_h(A) = \sum_{l = 1}^{\infty}\sum_{s \in C_l} h(s)\ind_{A}(s)
\]
for all $A \in \mathscr{P}(L)$.
\end{itemize}
\end{proposition}
We have arrived at a point now, where we can apply our general
framework to discontinuities of the first kind. In particular, we
can easily provide a representation of unordered sums over all jumps
of $f$; a fact, which frequently is used in the literature on
general semimartingales including references on L\'{e}vy processes,
but which seemingly hasn't been \textit{rigorously} proven yet, very
similar to the case of the proof of Theorem \ref{thm:at most
countably many jumps on compact intervals} (cf. e.\,g. \cite{Ap04},
\cite{Kl98}, \cite{Pa67}).
\begin{proposition}\label{prop:sum of jumps on R+ with invertible function}
Let $f : \R_+ \longrightarrow \R$ be an arbitrary function, and
assume that $\emptyset \not= J(f)$. Let $\{D_n : k \in \N\}$ be an
arbitrary partition of $J(f)$. Let $\Phi : \R_+ \longrightarrow
\R_+$ be strictly increasing and continuous such that $\Phi(0) = 0$.
Let $B$ be a non-empty subset of $L(f)$. Then
\[
\sum_{s \in B} \Phi\big(\vert \Delta f(s) \vert\big) = \mu_{\Phi
\circ \vert \Delta f \vert}\big(B \cap J(f) \big) = \sum_{n =
1}^{\infty} \sum_{s \in D_n} \Phi\big(\vert \Delta f(s) \vert\big)
\ind_{B}(s) = \sum_{n = 1}^{\infty} \Phi\big(\vert \Delta
f(\varphi(n)) \vert\big)
\]
for all bijective mappings $\varphi: \N \longrightarrow B \cap
J(f)$. If in addition $f$ is regulated, then
\[
\sum_{0 < s \leq t} \Phi\big(\vert \Delta f(s) \vert\big) =
\mu_{\Phi \circ \vert \Delta f \vert}\big((0, t] \cap J(f)\big) =
\sum_{n = 1}^{\infty} \sum_{s \in D_n} \Phi\big(\vert \Delta f(s)
\vert\big) \ind_{(a, t)}(s) = \sum_{n = 1}^{\infty} \Phi\big(\vert
\Delta f(\varphi(n)) \vert\big)
\]
for all $t \in (0, \infty)$ and bijective mappings $\varphi: \N
\longrightarrow (0, t] \cap J(f)$.
\end{proposition}
\begin{proof}
Let $B$ be an arbitrary non-empty subset of $L(f)$. Due to Theorem
\ref{thm:at most countably many jumps on R+}, $D : = J(f)$ is at
most countable. Since $\Phi : \R_+ \longrightarrow \R_+$ is strictly
increasing and continuous, it is invertible, and $\Phi^{-1} : \R_+
\longrightarrow \R_+$ is strictly increasing as well (due to the
Inverse Function Theorem). Hence, the set $\{s \in B: \Phi(\vert
\Delta f(s) \vert) > 0 \} = \{s \in B: \vert \Delta f(s) \vert > 0
\} = B \cap D$ is an at most countable subset of $D \subseteq L(f)$.
Consider the
the function $h : = \Phi \circ \vert \Delta f \vert$. Then
$\mu_h(B \cap D) \stackrel{\eqref{eq:mu of B cap D}}{=}
\mu_{h\iind{2}_B}(D).$
Since $\Phi(0) = 0$, equality \eqref{eq:mu of B cap D} implies that
$\mu_h(B \cap (L\setminus D)) = \mu_{h\iind{2}_B}(L\setminus D) =
0$, and it follows that
\[
\sum_{s \in B}h(s) = \mu_{h}(B) = \mu_{h}(B \cap D) =
\mu_{h\iind{2}_B}(D) = \sum_{n = 1}^{\infty}\sum_{s \in D_n} h(s)
\ind_{B}(s).
\]
Since the set $B \cap D$ is an at most countable subset of $L(f)$,
the first statement follows by Corollary \ref{cor:countable sums}.
If in addition $f$ is regulated, then $L(f) = (0, \infty)$ (due to
Remark \ref{remark: regulated fct}), implying that $B : = (0, t]
\subseteq L(f)$ for all $t \in (0, \infty)$. Now, the second
statement follows immediately from the first one. \qed
\end{proof}
If $f : \R_+ \longrightarrow \R$ were regulated, a natural question
would be to ask for the representation of the function of jumps
$\Delta g : L(g) \longrightarrow \R$, where $g(t) : = \sum_{s \in
(0, t]} \Phi\big(\vert \Delta f(s) \vert\big)$, $t \in (0, \infty)$.
To this end, let $h : (0,\infty) \longrightarrow \R_+$ be an
arbitrary positive real-valued function, and assume that $\mu_h((0,
\infty)) = \sum_{s \in (0, \infty)}h(s) < \infty$. Then $g(t) : =
\sum_{s \in (0, t]}h(s) = \mu_h((0, t]) \cap P) < +\infty$ for all
$t \in \R_+$, where $P : = \{s \in (0, \infty) : h(s) > 0\}$. Let $t
\in (0, \infty)$. Since $\mu_h : \mathscr{P}((0,\infty))
\longrightarrow \R_+$ is a (finite) measure, it follows that $g(t +
\frac{1}{n}) - g(t - \frac{1}{n}) = \mu_h(I_n)$ for sufficiently
large $n \in \N$, where $I_n : = (t - \frac{1}{n}, t +
\frac{1}{n}]$. Since $I_n \downarrow \{t\}$ as $n \to \infty$, we
obviously have
\[
\Delta g(t) = \lim_{n \to \infty} \mu_h(I_n) = \mu_h(\{t\}) = h(t).
\]
Hence, $L(g) = (0, \infty)$, and $\Delta \Big(\sum_{s \in (0, \cdot
)}h(s)\Big)= \Delta g = h$ on $(0, \infty)$. Moreover, since
$g(\frac{1}{n}) = \mu_h((0, \frac{1}{n}]) \to \mu_h(\emptyset)= 0$
as $n \to \infty$, it follows that $0 \in L^+(g)$ and $g(0+) = 0 =
g(0)$. Consequently, Remark \ref{remark: regulated fct}
implies the following
\begin{proposition}
Let $h : (0,\infty) \longrightarrow \R_+$ be an arbitrary positive
real-valued function. If $\sum_{s \in (0, \infty)}h(s) < \infty$,
then the function
\begin{eqnarray*}
g : \R_+ & \longrightarrow & \R_+\\
t & \mapsto & \sum_{s \in (0, t]} h(s),
\end{eqnarray*}
is regulated, $L(g) = (0, \infty)$, and
\[
\Delta g = h.
\]
\end{proposition}
\section{Random measures induced by optional processes}
Next, we transfer the main results of our previous investigations
to (trajectories of) stochastic processes. Again, let $(\Omega ,
\mathcal{F}, {\bf{F}}, {\Bbb{P}})$ be a given filtered probability
space such that ${\bf{F}}$ satisfies the usual conditions. If $X$
is an adapted and c\`{a}dl\`{a}g process, then it is well-known
that the left limit process $X_{-}$ is predictable. Recall that
every adapted and right-continuous process is optional (cf.
Theorem \ref{thm:POPA}). Consequently, we deal with a special case
of the slightly more general
\begin{lemma}\label{lemma:optional and regulated paths}
Let $(\Omega , \mathcal{F}, {\bf{F}}, {\Bbb{P}})$ be a filtered
probability space such that ${\bf{F}}$ satisfies the usual
conditions. Let $X : \Omega \times \R_+ \longrightarrow \R$ be a
stochastic process such that all trajectories of $X$ are regulated.
Then all trajectories of the left limit process $X_{-}$
$($respectively of the right limit process $X_{+}$$)$ are
left-continuous $($respectively right-continuous$)$. If in addition
$X$ is optional, then $X_{-}$ is predictable.
\end{lemma}
\begin{proof}
Fix $\omega \in \Omega$ and consider the (fixed) trajectory $f : =
X_{\bullet}(\omega)  : \R_+ \ \longrightarrow \R$. Since $f$ is a
regulated function, it follows that that $L^{-}(f) = (0, \infty)$.
Consequently, due to Corollary \ref{cor:f+ and f- for regulated f},
it clearly follows that the trajectory $Y_{\bullet}(\omega)$ of the
left limit process $Y : = X_{-}$ is left-continuous on $(0,
\infty)$. Similarly, it follows that the trajectory of the right
limit process $X_{+}$ is right-continuous on $\R_+ = L^{+}(f)$. Now
assume that in addition $X$ is optional and therefore adapted. Then
$X_{-}$ is an adapted process too. Consequently it follows that $Y =
X_{-}$ is adapted and left-continuous, and the definition of
predictability finishes the proof.\qed
\end{proof}
Now, we return to \textit{finitely} layered partitions and start
with the following observation. Despite its seemingly clear
context, it will be of high importance for our further
investigations.
\begin{lemma}\label{lemma:recursive construction of finite sets}
Let $\emptyset \not= D$ be a finite subset of $\R$, consisting of
$\kappa_D$ elements. Consider
\[
s_1^D : = \min(D)
\]
and, if $\kappa_D \geq 2$,
\[
s_{n+1}^D := \min(D \cap (s_n^D, \infty)\big),
\]
where $n \in \{1, 2, \ldots, \kappa_D - 1\}$. Then $D \cap
(s_{n}^D, \infty) \not= \emptyset$ for all $n \in \{1, 2, \ldots,
\kappa_D - 1\}$ and $s_{n}^D < s_{n+1}^D$ for all $n \in \{1, 2,
\ldots, \kappa_D - 1\}$. Moreover, we have
\[
D = \cdju{n=1}{\kappa_D}\{s_n^D\}.
\]
\end{lemma}
\begin{proof}
Obviously, nothing is to prove if $\kappa_D \in \{1, 2\}$. Let
$\kappa_D \geq 3$. Obviously, we have $D \cap (s_{1}^D, \infty)
\not= \emptyset$. Now assume by contradiction that there exists $n
\in \{2, \ldots, \kappa_D - 1\}$ such that $D \cap (s_{n}^D,
\infty) = \emptyset$. Choose the \textit{minimal} $m \in \{2,
\ldots, \kappa_D - 1\}$ such that $D \cap (s_{m}^D, \infty) =
\emptyset$. Then $s_{k}^D : = \min(D \cap (s_{k-1}^D, \infty)\big)
\in D$ is well-defined for all $k \in \{2, \ldots, m\}$, and we
obviously have $s_{1}^D < s_{2}^D < \ldots < s_{m}^D$. Moreover,
by construction of $m$, it follows that
\begin{equation}\label{eq:induction statement}
s \leq s_{m}^D \textup{ for all } s \in D.
\end{equation}
Assume now that there exists $s^\ast \in D$ such that $s^\ast
\not\in \{s_{1}^D, s_{2}^D, \ldots, s_{m}^D\}$. Then, by
\eqref{eq:induction statement}, there must exist $l \in \{1, 2,
\ldots, m -1 \}$ such that $s_{l}^D < s^\ast < s_{l+1}^D$, which is
a contradiction, due to the definition of $s_{l+1}^D$. Hence, such a
value $s^\ast$ cannot exist, and it consequently follows that $D =
\{s_{1}^D, s_{2}^D, \ldots, s_{m}^D\}$. But then $m =
\textup{card}(D) \leq \kappa_D - 1 < \kappa_D$, which is a
contradiction. Hence, $s_{k}^D \in D$ is well-defined for all $k \in
\{1, 2, \ldots, \kappa_D\}$. Since $\textup{card}(D) = \kappa_D$,
the proof is finished. \qed
\end{proof}
\begin{theorem}\label{thm: jumps of optional processes and stopping times}
Let $(\Omega , \mathcal{F}, {\bf{F}}, {\Bbb{P}})$ be a filtered
probability space such that ${\bf{F}}$ satisfies the usual
conditions. Let $X : \Omega \times \R_+ \longrightarrow \R$ be an
optional process such that all trajectories of $X$ are regulated.
Then $\Delta X$ is also optional. Put $X_{0-} : = X_{0+}$. If for
each trajectory of $X$ its set of jumps is not finite, then there
exists a sequence of stopping times $(T_n)_{n \in \N}$ such that
$(T_n(\omega))_{n \in \N}$ is a strictly increasing sequence in $(0,
\infty)$ for all $\omega \in \Omega$ and
\[
J(X_{\bullet}(\omega)) = \cdju{n=1}{\infty}\{T_n(\omega)\} \textup{
for all } \omega \in \Omega,
\]
or equivalently,
\[
\{\Delta X \not= 0\} = \cdju{n = 1}{\infty} \rsi T_n \lsi.
\]
\end{theorem}
\begin{proof}
Since the filtration is right-continuous, a direct calculation shows
that the right limit process $X_{+}$ is adapted. Due to Lemma
\ref{lemma:optional and regulated paths}, all paths of $X_{+}$ are
right-continuous on $\R_+$. Since the filtration is
right-continuous, it follows that $X_{+}$ is also adapted and hence
an optional process. Consequently, since the process $X$ was assumed
to be optional, and since each predictable process is optional, a
further application of Lemma \ref{lemma:optional and regulated
paths} implies that the jump process $\Delta X = X_{+} - X_{-}$ is
the sum of two optional processes, hence optional itself.
Fix $\omega \in \Omega$. Consider the trajectory $f : =
X_{\bullet}(\omega)$. Due to statement (ii) of Theorem \ref{thm:at
most countably many jumps on R+}, there exists a finitely layered
partition of $J(f) \subseteq L(f) = (0, \infty)$ which now is
randomised, and it follows that we may write $J(f)$ as
\[
J(f) = \dju{m = 1}{\infty}{4.6}{2.2} D_m(\omega),
\]
where $\kappa_m(\omega) : = \textup{card}(D_m(\omega)) < +\infty$
for all $m \in \N$.
Let ${\Bbb{M(\omega)}} : =\{m \in \N : D_{m}(\omega) \not=
\emptyset\}$. Fix an arbitrary $m \in {\Bbb{M(\omega)}}$. Consider
\[
0 < S_1^{(m)}(\omega) : = \min(D_{m}(\omega))
\]
and, if $\kappa_m(\omega) \geq 2$,
\[
0 < S_{n+1}^{(m)}(\omega) := \min\big(D_{m}(\omega) \cap
(S_n^{(m)}(\omega), \infty)\big),
\]
where $n \in \{1, 2, \ldots, \kappa_m(\omega) - 1\}$.
Since $\Delta X$ is optional, it follows that $\{\Delta X \in B\}$
is optional for all Borel sets $B \in \mathcal{B}(\R)$. Moreover,
since $\Delta f(0) = \Delta X_0(\omega) : = 0$ (by assumption), it
actually follows that $\{s \in \R_+ : (\omega, s) \in \{\Delta X \in
C\}\} = \{s \in (0, \infty) : (\omega, s) \in \{\Delta X \in C\}\}$
for all Borel sets $C \in \mathcal{B}(\R)$ which do not contain $0$.
Hence, as the construction of the sets $D_{m}(\omega)$ in the proof
of Theorem \ref{thm:at most countably many jumps on R+} clearly
shows, $S_1^{(m)}$ is the d\'{e}but of an optional set.
Consequently, due to Theorem \ref{thm:optional debut is stopping
time}, it follows that $S_1^{(m)}$ is a stopping time. If
$S_n^{(m)}$ is a stopping time, the stochastic interval $\lsi
S_n^{(m)}, \infty \rsi$ is optional too (cf. \cite{HeWaYa92},
Theorem 3.16). Thus, by construction, $S_{n+1}^{(m)}$ is the
d\'{e}but of an optional set and hence a stopping time. Due to Lemma
\ref{lemma:recursive construction of finite sets}, we have
\begin{equation}\label{eq:rep of the set of jumps of an optional regulated
process} J(f) = \dju{m \in {\Bbb{M(\omega)}}}{}{7}{2.2}
D_{m}(\omega) = \dju{m \in {\Bbb{M(\omega)}}}{}{7}{4.2} \dju{n =
1}{\kappa_m(\omega)}{5.5}{2.2}\{S_n^{(m)}(\omega)\}.
\end{equation}
Hence, since for each trajectory of $X$ its set of jumps is not
finite, the at most countable set ${\Bbb{M(\omega)}}$ is not finite,
hence countable, and a simple relabelling of the stopping times
$S_n^{(m)}$ finishes the proof. \qed
\end{proof}
We will see now how the choice of finitely layered jump partitions
enables a natural approach for recovering the jump measure of a
c\`{a}dl\`{a}g and adapted stochastic process, and how we can
transfer the structure of the jump measure to more general classes
of optional processes which need not necessarily to be
right-continuous.
\begin{proposition}\label{prop:random measure for general processes X}
Let $(\Omega , \mathcal{F}, {\bf{F}}, {\Bbb{P}})$ be a filtered
probability space such that ${\bf{F}}$ satisfies the usual
conditions. Let $X : \Omega \times \R_+ \longrightarrow \R$ be an
optional process such that all trajectories of $X$ are regulated and
$\Delta X_0 : = 0$. Consider
\[
M_X(\omega) : = \Big\{(s, \Delta X_{s}(\omega)) : s \in
J(X_{\bullet}(\omega)) \Big\},
\]
where $\omega \in \Omega$. Then
\[
\textup{card}(M_X(\omega) \cap G) = \sum_{s \in
J(X_{\bullet}(\omega))} \ind_{G}\big(s, \Delta X_{s}(\omega) \big) =
\sum_{s > 0} \ind_{G}\big(s, \Delta X_{s}(\omega) \big)
\ind_{\{\Delta X \not= 0\}}(\omega, s)
\]
for all $(\omega, G) \in \Omega \times \mathcal{B}(\R_+) \otimes
\mathcal{B}(\R)$.
\end{proposition}
\begin{proof} For simplicity reasons, we may assume that the set of
jumps of each trajectory of $X$ is not finite (due to Theorem
\ref{thm:Fubini for unordered sums} and representation \eqref{eq:rep
of the set of jumps of an optional regulated process}). Fix
$(\omega, G) \in \Omega \times \mathcal{B}(\R_+) \otimes
\mathcal{B}(\R)$. Consider $D_X(\omega) : = J(X_{\bullet}(\omega))
\subseteq L(X_{\bullet}(\omega)) = (0, \infty)$. Theorem \ref{thm:at
most countably many jumps on R+} implies that the
\textcolor{ownred}{random set}
\begin{eqnarray*}
M_X(\omega) : = \Big\{(s, \Delta X_{s}(\omega)) : s \in
D_X(\omega) \Big\} = \mbox{Gr}\Big(\Delta
X_{\bullet}(\omega)\vert_{D_X(\omega)}\Big)
\end{eqnarray*}
\textcolor{owngreen}{is at most countable}.
Put $j\deep{X}(\omega, \textcolor{owngreen2}{G}) : =
\textup{card}(M_X(\omega) \cap G)$. Due to Theorem \ref{thm: jumps
of optional processes and stopping times}, it follows that there
exists a sequence of stopping times $\{T_n : n \in \N\}$ such that
\[ M_X(\omega) \cap \textcolor{owngreen2}{G} = \dju{m =
1}{\infty}{4.5}{2.2} \big\{\big(T_m(\omega), \Delta X_{T_m}(\omega)
\big)\big\} \cap G,
\]
where $\Delta X_{T_m}(\omega) : = \Delta X_{T_m(\omega)}(\omega) =
\Delta X(\omega, T_m(\omega))$. Since all unions are disjoint ones,
the respective cardinals are additive. Consequently,
\[
j\deep{X}(\omega, \textcolor{owngreen2}{G}) = \sum_{m = 1}^\infty
\mbox{card}\big( \big\{\big(T_m(\omega), \Delta X_{T_m}(\omega)
\big)\big\} \cap G \big) = \sum_{m = 1}^\infty
\ind_{G}\big(T_m(\omega), \Delta X_{T_m}(\omega) \big),
\]
and Theorem \ref{thm:Fubini for unordered sums} together with Lemma
\ref{lemma:indicator function in sum} imply that
\[
j\deep{X}(\omega, \textcolor{owngreen2}{G}) = \sum_{s \in
D_X(\omega)} \ind_{G}\big(s, \Delta X_{s}(\omega) \big) = \sum_{s >
0} \ind_{G}\big(s, \Delta X_{s}(\omega) \big) \ind_{\{\Delta X \not=
0\}}(\omega, s).
\]
\qed
\end{proof}
\begin{theorem}\label{thm:optional random measures}
Let $(\Omega , \mathcal{F}, {\bf{F}}, {\Bbb{P}})$ be a filtered
probability space such that ${\bf{F}}$ satisfies the usual
conditions. Let $X : \Omega \times \R_+ \longrightarrow \R$ be an
optional process such that all trajectories of $X$ are regulated and
$\Delta X_0 : = 0$. Then the function
\begin{eqnarray*}
j\deep{X} : \Omega \times \mathcal{B}(\R_+) \otimes
\mathcal{B}(\R) & \longrightarrow & \Z_+ \cup \{+ \infty\}\\
(\omega, G) & \mapsto & \sum_{s > 0} \ind_{G}\big(s, \Delta
X_{s}(\omega) \big) \ind_{\{\Delta X \not= 0\}}(\omega, s)
\end{eqnarray*}
is an integer-valued random measure.
\end{theorem}
\begin{proof}
We only have to combine Theorem \ref{thm: jumps of optional
processes and stopping times} and \cite{HeWaYa92}, Theorem
11.13.\qed
\end{proof}
Put $G : = B \times \Lambda$, where $B \in \mathcal{B}(\R_+)$ and
$\Lambda \in \mathcal{B}(\R)$. Since each trajectory of $X$ is
regulated, $B\setminus\{0\} \subseteq (0, \infty) =
L(X_\bullet(\omega))$. Hence, Lemma \ref{lemma:indicator function in
sum} directly leads to the following representation:
\[
j\deep{X}(\omega, B \times \Lambda) = j\deep{X}(\omega,
(B\setminus\{0\}) \times \Lambda) = \sum_{s \in B\setminus\{0\}}
\ind_\Lambda(\Delta X_s(\omega))\ind_{\{\Delta X \not= 0\}}(\omega,
s)
\]
for all $B \in \mathcal{B}(\R_+)$ and $\Lambda \in
\mathcal{B}(\R)$.

To sum up, $j\deep{X}(\omega, G)$ in general counts,
$\omega$-by-$\omega$, the number of all $s > 0$ such that $\Delta
X\deep{s}(\omega) \not= 0$ and $(s, \Delta X_s(\omega)) \in G$. In
other words,
\[
j\deep{X}(\omega, (\mbox{d}t, \mbox{d}x)) = \sum_{s
>0}\ind_{\{\Delta X \not= 0\}}(\omega, s)\cdot\delta_{\big(s, \Delta
X_{s}(\omega)\big)}(\mbox{d}t,\mbox{d}x).
\]
We finish this paper by considering right-continuous trajectories
again and note the following
\begin{corollary}
Let $(\Omega , \mathcal{F}, {\bf{F}}, {\Bbb{P}})$ be a filtered
probability space such that ${\bf{F}}$ satisfies the usual
conditions. If $X : \Omega \times \R_+ \longrightarrow \R$ is an
adapted and c\`adl\`ag process such that $\Delta X_0 : = 0$, then
$j_X$ is the jump measure of $X$.
\end{corollary}

\textbf{Acknowledgements}: The author gratefully thanks Dave
Applebaum and Finbarr Holland for highly fruitful discussions
including Finbarr Holland's information about the very useful M. Sc.
Thesis \cite{Od04}.

\end{document}